\documentclass[preprint, 12pt]{elsarticle}
\usepackage{graphicx}
\usepackage{amsmath}
\usepackage{amssymb}

\newcommand{\pr}{\mbox{\sf P}}
\newcommand{\ex}{{\bf\sf E}}               

\newcommand{\bs}{{\bf s}}               

\newcommand{\bd}{{\bf d}}
\newcommand{\bbe}{{\bf e}}

\newtheorem{thm}{Theorem}

\newtheorem{pro}[thm]{Proposition}

\def\th{\theta}

\numberwithin{equation}{section}

\journal{Operations Research Letters}
\begin{document}

\title{A Stochastic Knapsack Game: Revenue Management in Competitions
}


\author{Yingdong Lu 
}


\date{Received: date / Accepted: date}

\begin{abstract}
We study a mathematical model for revenue management under competition with multiple sellers. The model combines the stochastic knapsack problem, a classic revenue management model, with a non coorperative game model that characterizes the sellers' rational behavior. We are able to establish a dynamic recursive procedure that incorporate the value function with the utility function of the games. The formalization of the dynamic recursion allows us to establish some fundamental structural properties.
\end{abstract}

\maketitle

\section{Introduction}
\label{sec:intro}

A key model in revenue (yield) management is the following, a seller needs to sell a fixed amount of certain commodity before a fixed deadline to different buyers with individual price they are willing to pay, and the seller can dynamically adjust the selling price to maximize his/her overall revenue over time.  Stochastic knapsack problem, also known as stochastic dynamic knapsack problem, a mathematical problem that captures the essences of this model,  quantifies some of the most fundamental trade-offs in revenue management, and serves as an important building block for more complex and sophisticated models for real life applications. Consequently, the stochastic (dynamic) knapsack problem and its variations have been studied extensively, see, e.g.~\cite{GallegovanRyzin1994},~\cite{BitranMondschein95},~\cite{FengGallego1995},~\cite{PapastavrouRajagopalanKleywegt1996},~\cite{FengGallego2000},~\cite{FengXiao2000},~\cite{VanSlykeYoung2000},~\cite{ZhaoZheng2000},~\cite{LinLuYao2008}. It is one of the fundamental models surveyed by Anould de Boer in~\cite{DenBoer2015}, please refer to that paper for more details, as well as references.

It is natural to ask the question of what would happen if there are multiple sellers competing for the same demand stream from the buyers. In this paper, we generalize the classic stochastic knapsack problem, and formulate a mathematical model to capture the basic relations in this situation. An immediate goal is to formulate a dynamic recursion for calculating the optimal policies for sellers. In the single seller case, this is accomplished through the formulation of a dynamic program that computes the maximum expected revenue starting at any time with any amount of remaining inventory. However, in the case of multiple sellers, at each time period, the sellers' decisions are inter-dependent. It is, therefore, not a trivial task to decide what will be the next best action even if every seller has the same forecast of the future demand arrivals. Another difficulty is that when multiple sellers are willing to sell the product, the buyer can have different ways to choose one of them to fulfill the demand, the difference in these selection rules has significant impact on the evolution of the system. To overcome these difficulties, we model the sellers as rational individual or institutions, and introduce a noncooperative game at each step of the dynamic recursion characterizing their behavior. Furthermore, we follow a static probabilistic selection rule, which will be described precisely later, that the buyer will use to select sellers. This selection rule, on one hand, reflects market power of the sellers, on the other hand, it allows the uncertainty that is natural in business reality.  With this mechanism, the utility functions of the games are properly connected with the value functions of the  dynamic recursion, thus help to identify pure strategy Nash equilibriums. Under the selection rule assumed, we are able to demonstrate that there is a unique Nash equilibrium of the game. In turn, assuming that the Nash equilibriums will be the strategy followed by all the sellers at each step, the dynamic recursion is able to proceed. Once establishing the dynamic recursion, we are able to extend the arguments that are effective for the single seller dynamic programming, and demonstrate that, in some cases, the value function exhibits remarkable rich monotonicity properties that provide insights to key trade-offs to the problem and can be helpful to dynamic pricing in practice.  A related but different model is considered in \cite{doi:10.1287/mnsc.2013.1821}, it is concluded that, under a differential game setting, the equilibrium structure enjoys simple structural properties. While the model studied here is quite different, but results are similar in spirit.

The rest of the paper will be organized as follow. In Sec. \ref{sec:models}, we will introduce the basic mathematical models, and review preliminaries including some basic concepts in game theorety that will be needed for our analysis. In Sec. \ref{sec:dp}, we will discuss in details the dynamic recursion in which the game aspect of the problem is incorporated. In Sec. \ref{sec:monotonicity}, we establish some fundamental structural properties of the value functions of the dynamic recursion. Finally, we conclude the paper in Sec. \ref{sec:conclusions} with a summary of our findings.

\section{Models and Preliminaries}
\label{sec:models}

\subsection{Model Descriptions}
\label{sec:model_descriptions}

Suppose that there are $N$ sellers, and each seller $n$,  $n=1,2,\ldots, N$, has an initial inventory of  $C_n$ units of product ( could be either goods or services) at the beginning of a common selling horizon. The selling horizon is discrete and of length $T<\infty$. At each time $t=1,2,\ldots, T$, demand for one unit of the product will emerge, and the buyer will post a price that he/she is willing to pay. To accommodate the event of no arrival, we can always include a class of demand with exceedingly low price. The sellers who have positive inventory need to decide whether they should accept or reject this demand. The buyer will then select one seller among all the sellers that accept the demand according to certain {\it selection rule}, and the selected seller will supply the product and collect the revenue. At the end of the selling horizon, all the remaining product will be savaged. The goal for each seller is to maximize his/her expected revenue. 

We assume that each seller does not have the information of the exact value of the initial inventory of other sellers, but has a distributional estimation of that quantity. We also assume that the distributional information of the future demand price is given to each seller, and no seller has any extra knowledge. In particular, we assume that the price of the demand realization at each period follow an independent and identically distributed discrete probability distribution $P$, with $\pr[P=p_i]= \th_i$, $i=1,2,\ldots,I$.

Suppose that, at each time $t$, when the demand is of class $i$,  i,e. the price is $p_i$, a subset of sellers, denoted by $A_t(i)$ (which can be shortened to $A_t$ when there is no ambiguity), will accept the demand, decided based on the remaining time, demand type, remaining inventory and the selling history up to time $t$. The buyer will select only one seller among them, which means that there is a possibility that no seller is selected. There could be various selection rule models reflecting different market mechanisms, for example, a static rule ( the buyer chooses one product over the other overwhelmingly, which happens often in some local and monopoly market) and weighted rule (buyer assigns weights to the each product, then randomly, with probabilities determined by the weights, select ones that are available). In this paper, we will focus on a random allocation rule with static probabilities: each seller is associated with a probability $\pi_n$, $\sum_{n=1}^N \pi_n=1$.  At each time, if a seller accepts, the probability of it being selected is always $\pi_n$, and with probability $1-\sum_{n \in A_t(i)}\pi_n$, no one is selected.

At each time $t$, the phenomenon that the sellers are making independent decisions based on distributional information on the other sellers can be best modeled by a non coorperative strategic game, see, e.g. \cite{osborne1994course}.

\section{A Dynamic Recursion Formulation}
\label{sec:dp}

Our goal is to identify a strategy for a seller to achieve the best outcome, in terms of average revenue, under a reasonable assumption on other sellers' behavior. Recall that each seller  $n$, $n=1,2,\ldots, N$ with initial inventory $C_n$ is also given the distributional information of the inventory of all other sellers, either though statistical forecast or other business information inquiry, and any two sellers will be given the exact same distribution on the third seller. In addition, all the sellers do observe all the sells outcomes up to each decision time epoch, i.e. they know the amount each seller sold so far. It is our intention to derive a dynamic recursion for calculating the best outcome, hence the optimal strategy for each seller. Equivalently, given $\bs=(s_1, s_2, \ldots, s_N)$ representing the amount of inventory has been sold so far by each seller, we seek to calculate $v_n(t, d_n, \bs)$,  $n=1,2,\ldots, N$, the maximum expected revenue seller $n$ can collect starting from time $t$ and with remaining inventory $d_n$, for any time $t=1,2,\ldots, T$.

We assume that the behavior of the sellers is modeled as a $N$-person game, and if sellers follow the Nash Equilibriums at each time period, a dynamic recursion can proceed. This will be argued inductively. At the last time period $T$, given a price realization, $p_i$, there are two strategies for each seller, accept or reject. The utility function of the game for seller $n$ will be the expected revenue collected by taking either action. If reject, of course, there is no revenue. 
It is clear that, if the random selection rule with static probabilities is followed, there exists a unique Nash equilibrium, that is, every seller will accept, as long as they have a positive inventory.  In this case, the value function $v_n(T, \bd)$ has the following form,
$v_n(T, \bd) = \pi_n \ex[P],$
where $\pi_n$  are the probabilities in the section rule model.

Now, suppose that we can calculate recursively all the value function $v_n(t+1, d_n, \bs)$ for any feasible $\bs$, we demonstrate that there exists a unique pure strategy Nash equilibrium at time period $t$,  and show how it is related to the calculation of the value function for time $t$, $v_n(t, d_n, \bs)$.
There are two actions for each seller, accept or reject. The payoff function will be the expected revenue to be collected until time $T$. Therefore, the seller will consider the following {\it balance inequality}, whose left hand side (LHS) represents the price we get immediately, and right hand side (RHS) represents the future reward,
\begin{align}
\label{eqn:balance}
p \ge \ex_{n, t}[ v_n (t+1, d_n, \bs)- v_n (t+1, d_n-1, \bs+\bbe_n)],
\end{align}
where $\ex_{n, t}$ is the expectation with respect to the information available at time $t$ for seller $n$, $p$ the generic price the class indicator is suppressed when there is no ambiguity). If \eqref{eqn:balance} holds, then the order will be accepted. Otherwise, if we have,
\begin{align}
\label{eqn:unbalance}
p < \ex_{n, t}[ v_n (t+1, d_n, \bs)- v_n (t+1, d_n-1, \bs+\bbe_n)],
\end{align}
the order will be rejected.

\noindent 
{\bf Remark}  The operator $\ex_{n,t}$ can be treated in a way as an conditional expectation, the information update each time is basically the confirmation that the random variable of each seller's inventory is larger than the cumulative sales, which is updated at the end of each time period. 

\begin{pro}
The above defined strategy is a unique Nash Equilibrium.
\end{pro}
{\bf Proof }To prove that it is a Nash Equilibrium, let us discuss separately for those sellers depends upon their decisions. Suppose that for a particular seller $n$, the action is to accept, three events can happen,
\begin{itemize}
	\item 
	sell $n$ is selected, with probability $\pi_n$;
	\item
	some other seller $j$ in $A_t$ is selected, with probability $\pi_j$;
	\item
	no seller is selected, with probability $1- \sum_{A_t} \pi_i$.
\end{itemize}
Sum them up, the pay-off function has the following form, 
\begin{align*}
	\pi_n\ex_{n, t}[p+v_{n} (t+1, d_n-1, \bs)]& + \sum_{j\neq n, j\in A_t} \pi_j \ex_{n, t}[v_n(t+1, d_n, \bs +\bbe_j)]\\ &+ \left(1- \sum_{A_t} \pi_i\right) \ex_{n, t}[v_n(t+1, d_n, \bs)].
\end{align*}
If seller $n$ deviates from the strategy, i.e., rejects the demand, its payoff will be, $$\sum_{j\neq n, j\in A_t} \pi_j \ex_{n, t}[v_n(t+1, d_n, \bs +\bbe_j)]+ \left(1- \sum_{A_t-\{n\}} \pi_i\right) \ex_{n, t}[v_n(t+1, d_n, \bs)].$$ 
From the \eqref{eqn:balance}, we know that seller $n$ could not be better off. 

In the case seller $n$ reject, the pay off is, 
\begin{align*}\sum_{j\neq n, j\in A_t} \pi_j \ex_{n, t} [v_n(t+1, d_n, \bs +\bbe_j)]+ \left(1- \sum_{A_t-\{n\}} \pi_i\right) \ex_{n, t}[v_n(t+1, d_n, \bs +\bbe_j)].\end{align*}
If the seller deviates from this strategy, the pay-off will become,
\begin{align*}
\pi_n[p+ \ex_{n, t}[v_{n} (t+1, d_n-1, \bs)]] &+ \sum_{j\neq n, j\in A_t} \pi_j \ex_{n, t} [v_n(t+1, d_n, \bs +\bbe_j)]\\ &+ \left(1- \sum_{A_t} \pi_i\right) \ex_{n, t} [v_n(t+1, d_n, \bs )].
\end{align*}
However, we know that $p+\ex_{n, t} [v_n(t+1, d_n-1, \bs )] <\ex_{n, t} [v_n(t+1, d_n, \bs)] $, therefore, the seller will be worse off. 

Suppose any other strategy that has a seller $n$, such that, $p+v_n(t+1, d_n-1, \bs +\bbe_n) <v_n(t+1, d_n \bs) $,  but seller $n$ accepts the demand. We can see that deviation will lead to better pay-off. Meanwhile if there is a seller $n$ with  $p+v_n(t+1, d_n-1, \bs +\bbe_n) \ge v_n(t+1, d_n, \bs)$, but seller $n$ rejects, a deviation will lead to higher pay-off.
$\Box$

The above arguments allow us to present the following dynamic recursion for the value function, 
\begin{align}
&v_n(t, d_n, \bs) =  \sum_{i=1}^I \th_iw_n(t+1, d_n,\bs, p_i), \label{eqn:main_recursion_1} \\
&w_n(t+1,  d_n,\bs, p_i)  = \left(1-\sum_{m=1}^N \pi_m\right)\ex_{n,t}[v_n(t+1,d_n, \bs)]\nonumber \\ &+ p_i \pi_n{\bf 1}\left\{\ex_{n,t}[v_n(t+1, d_n-1, \bs+\bbe_n)+p_i] \ge \ex_{n,t}[v_n(t+1, d_n, \bs)\right]\} \nonumber \\ & +\sum_{m=1}^N \pi_m [v_n(t+1, d_n, \bs +\bbe_m){\bf 1}\{\ex_{m,t}[v_m(t+1, d_m -1, \bs+\bbe_m)+p_i] \ge \ex_{m, t}[v_m(t+1, d_m)]\}  \nonumber \\ &  +\ex_{n,t}[v_n(t+1, d_n+\bs ){\bf 1}\{\ex_{n,t}[v_n(t+1, d_n-1, \bs+\bbe_n)+p_i < v_n(t+1, d_n, \bs)]\} ],\label{eqn:main_recursion_2}
  \\
 & v_N(T, d_n, \bs) =  \pi_N \ex[p].\label{eqn:main_recursion_3}
\end{align}

\noindent 
    {\bf Remark }
The information available at time $t$ is on the distribution on the initial capacity of all the other sellers, as well as the sales records in the past period. At time $t+1$, the sales records will be amended with what happened during time period $t$, the distribution inform hence is naturally updated, for example, if the original distributional estimation is $D$, and at time $t$, the total sales has been $s$, then that information should be updated to $D;D\ge s$. At time $t$, if there are sales by that seller, it should be updated to $D;D\ge s+1$, otherwise, it will stay at $D;D\ge s$.

\section{Monotonicity of the Value Functions and its Implications in Revenue Management}
\label{sec:monotonicity}

In this section, we will establish monotonicity properties of the value function $v_n(t, d, \bs)$, based on the dynamic recursion formulated in Sec. \ref{sec:dp}.
The main result is stated in the following theorem, and its proof is presented in the Appendix. 
\begin{thm}
\label{thm:monotone} Under the random selection rule with static probabilities, the value function of the knapsack problem $v_n(t, d, \bs)$ for the $n$-th seller satisfies the following monotonicity properties. 
\begin{itemize}
\item[(1)] Monotone in inventory $\bd$, i.e. $\ex_{n, t-1}[v_n(t, d_n, \bs) \ge \ex_{n, t-1}[v_n(t, d_n-1, \bs)]$;
\item[(2)] Monotone in selling amount of competitors,$$\ex_{n, t-1}[v_n(t, d_n, \bs)] \le \ex_{n, t-1}[v_n(t, d_n, \bs+\bbe_i)] ;$$
\item[(3)] Monotone in time $t$, i.e.  $\ex_{n, t-1}[v_n(t, d_n,\bs)] \ge \ex_{n, t}[v_n(t+1, d_n, \bs)]$
\item [(4)]"Concave" in $d_n$, i.e.,
\begin{align}&\ex_{n, t-1}[v_n(t,d, \bs)-v_n(t,d_n-1,\bs+\bbe_n)]\nonumber \\ \ge & \ex_{n, t-1}[v_n(t,d_n+1,\bs)-v_n(t,d, \bs+\bbe_n)];\label{eqn:concavity}\end{align}
\item[(5)]
Submodular in $(t,\bd)$, i.e.,
\begin{align}\label{eqn:submodularity}
&\ex_{n, t-1}[v_n(t,d, \bs)]-\ex_{n, t-1}[v_n(t,d, \bs+{\bf e}_n)]\nonumber \\ & \ge \ex_{n, t}[v_n(t+1,d, \bs)]-\ex_{n, t}[v_n(t+1,d,\bs+{\bf e}_n)].\end{align}  
\item[(6)] Submodular in $\bd$, i.e.
\begin{align}\label{eqn:submodularity_more}
&\ex_{n, t-1}[v_n(t,d,\bs)]-\ex_{n, t-1}[v_n(t,d, \bs+{\bf e}_n)]\nonumber \\ & \ge \ex_{n, t-1}[v_n(t,d,\bs-{\bf e}_m)]+\ex_{n, t-1}[v_n(t,d,\bs+{\bf e}_n-{\bf e}_m)],\end{align}
with $m\neq n$.
\end{itemize} 
\end{thm}


Recall that we raised several questions in the introduction, here, after the statement of the main monotonicity results, we need to use them to answer some of those questions. 

From (2) of Theorem \ref{thm:monotone}, we can immediately see that, 
\begin{pro}
For each individual seller, his/her total average revenue is a monotone decrease function of his/her competitors inventory surplus levels. 
\end{pro}
{\bf Remark }
It is apparent that the more the overall supply is, the less is the expected marginal gain for each individual unit.

From (4), i.e. "concave in inventory" of Theorem \ref{thm:monotone}, we can conclude that
\begin{pro}
If it is optimal to accept the at certain point, then it is also optimal to accept when your competitors have more inventory.
\end{pro}
{\bf  Remark}
The intuition is that when there are more inventory in the hands of the competitors, they will be more aggressive, and it will then lower your expected marginal gain. Thus, you will be more likely to accept a lower price.

The inequality in (6) tells us that 
\begin{pro}
A lower selling amount of his/her competitors will make a seller less likely to accept a fixed price; certainly, a higher selling amount will make the same seller more likely to accept the same price.  
\end{pro}
Intuitively, observing more sells from ones competitors will make a seller more aggressive.

\noindent

\section{Conclusions}
\label{sec:conclusions}

In this paper, we extend the classic stochastic knapsack problem to model competitions between several sellers and effects on their dynamic pricing decisions. By utilizing dynamic programming techniques, together with a game theoretical model on the sellers' behavior, we are able to identify a simple strategy, i.e. checking the balance inequality, for each seller, and a dynamic recursion for calculating the value functions required. Furthermore, we show that the value functions have several important first and second order monotonicity properties that are of important theoretical values and critical practical implications.

\begin{appendix}

\section{Proof of Theorem \ref{thm:monotone} }
\label{sec:proof}


\noindent
{\bf Proof}
It is easy to see that (1) and (3) are trivial. We will prove the rest by backward induction on time $t$. First, it is trivial to check all of them at the end of selling season, time $T$. Next, suppose that at time period $t+1$ and later, the properties (2) and(4) through (6) hold. We want to extend all the result to time period $t$. 
Since selection is based on the static probabilities, to facilitate our discussion, denote $\Pi_0$ the event that no seller is selected, and $\Pi_i, i=1,\ldots, N$ the event that seller $i$ is selected. From our assumptions, it is clear that the probabilities of these events are $\pi_i, i=0,1,\ldots, N$, respectively. Furthermore, since all the demand random variables are i.i.d, it is suffice to focus on the event that the price of the demand is $p_i$, $i=1,\ldots, N$. We will use a generic notation $p$ to denote the price, for the ease of exposition. 

\subsection*{Validity of (2)}

Recall that, we need to establish $\ex_{i, t-1}[v_i(t, d_i, \bs)] \le \ex_{i, t-1}[v_i(t, d_i, \bs+\bbe_j)]$,  for $j \neq i$. Without loss of generality, it suffices to show, $\ex_{1, t-1}[v_1(t, d_1, \bs) ]\le \ex_{1, t-1}[v_1(t, d_1, \bs+\bbe_j)]$, for any $j>1$. We will argue that the inequality holds on each event $\Pi_i, i=0,1,\ldots, N$.
On $\Pi_0$, since no seller is selected, it is easy to see that the inequality holds by induction, and the induction arguments also applies to $\Pi_i, i\neq 1$ and $i\neq j$. On $\Pi_1$, examine what happens at time $t$, the only case that is not straightforward is that seller one only accept given that the history is $\bs$ but reject when it is $\bs+\bbe_j$. In this case, the left hand side (LHS) of the inequality becomes $\ex_{1, t}[v_1(t+1, d_1-1, \bs+\bbe_1)]+p$.  By induction, it is less than or equal to $\ex_{1, t}[v_1(t+1, d_1-1, \bs+\bbe_1+\bbe_j)]+p$. Meanwhile, $\ex_{1, t}[v_1(t+1, d_1-1, \bs+\bbe_j+\bbe_1))]+p\le \ex_{1, t}[v_1(t+1, d_1, \bs+\bbe_j)]$ due to the fact that this demand is not accepted when the history is $\bs+\bbe_j$. Hence, the inequality follows. On $\Pi_j, j>1$, there are two cases need to be considered depending on whether seller $j$ accepts the demand. Case I,  seller $j$ only accepts when the history is $\bs$ not when it is $\bs+\bbe_j$. In this case, we have both the LHS and the right hand side (RHS) equal to $\ex_{1, t}[v_1(t+1, d_1, \bs+\bbe_j)]$.
Case II, seller $j$ accepts in both cases. Then, the desired inequality is a consequence of $\ex_{1, t}[v_1(t+1, d_1, \bs+\bbe_1)] \le \ex_{1, t}[v_1(t+1, d_1, \bs+2\bbe_1)]$, which is the consequence of induction.

\subsection*{Validity of (4)}

Without loss of generality, we only need to show, 
\begin{align*}
&\ex_{1, t-1}[v_1(t, d_1, \bs)] -\ex_{1, t-1}[v_1 (t, d_1 -1, \bs+\bbe_1)]\\ \ge &\ex_{1, t-1}[v_1(t, d_1+1, \bs )] - \ex_{1, t-1}[v_1(t, d_1, \bs+\bbe_1)]. 
\end{align*}
Let us first consider case by case based on whether demand will be accepted by seller one. From the induction assumption for time $t+1$, we know that there are only the following cases, 
\begin{itemize}
\item[I.]
the demand is only accepted when the inventory is at $d_1+1$ not when it is $d_1$;
\item[II.]
the demand is accepted when the inventory levels are at  both $d_1+1$ and $ d_1$;
\item[III.]
the demand is rejected in either case. 
\end{itemize}
And we will discuss each case for events $\Pi_0, \Pi_1$ and $\Pi_j, j>1$.

In Case I, on event $\Pi_0$, the inequality follows from induction, i.e. the concavity with respect to the inventory, at time $t+1$. On the event $\Pi_1$, the LHS becomes $\ex_{1, t}[v_1(t+1, d_1, \bs) -v_1 (t+1, d_1-1, \bs+\bbe_1)] $, and the RHS becomes $p$, then the inequality follows because the balance inequality is violated, which is exactly the reason the demand is not accepted when the inventory is at $(d_1, \bs)$.  On $\Pi_j$, $j \ge 2$, since the decision of seller $j$ will not depend on the actual amount of inventory seller one has, but just the distribution, the RHS becomes, $ \ex_{1, t}[v_1(t+1, d_1+1, \bs+\bbe_j ) - v_1(t+1, d_1, \bs+\bbe_j+\bbe_1)]$. Hence, the inequality will follow from the concavity with respect to inventory from time $t+1$ due to induction assumption. 

In Case II, again, we only need to look at event $\Pi_1$, where the LHS becomes $p$ and the RHS becomes $\ex_{1, t}[v_1(t+1, d_1, \bs) -v_1 (t+1, d_1-1, \bs+\bbe_1)]$, and the inequality follows from the balance inequality. Finally, in Case III, the inequality follows from induction.


\subsection*{Validity of (5) }
Again, we need to show that,
\begin{align*}
&\ex_{1, t-1}[v_1(t, d_1, \bs)] - \ex_{1, t-1}[v_1(t, d_1-1, \bs+\bbe_1)] \\ \ge &  \ex_{1, t-1}[v_1(t+1, d_1, \bs)] - \ex_{1, t-1}[v_1(t+1, d_1-1, \bs+\bbe_1)].
\end{align*}
We will examine the inequality on each event $\Pi_i$, i=0,1, \ldots, N.
On $\Pi_0$, the inequality follows directly from the induction assumption. 
On $\Pi_1$, let us consider three subcases. First, it is again a straightforward conclusion from the induction assumption if the demand is not accepted for either inventory level. On the other hand if it is accepted for both inventory levels, then the inequality holds due to the induction assumption on the validity of (4) at time $t$ and $t+1$. If seller one only accepts when the inventory level is at $d_1$, but not when it is at $d_1-1$, the LHS will become $p$, then by the condition of accept, i.e. the balance inequality, it is larger than the RHS.
On $\Pi_j, j \ge 2$, the inequality follows from the induction assumption on (6) if the demand is accepted for both inventory levels. By the distributional assumption, that is all that needs to be considered. 

\subsection*{Validity of (6) }
It is  our task to show that, for $j\ge 2$, 
\begin{align*}
&\ex_{1, t-1}[v_1(t, d_1, \bs)] - \ex_{1, t-1}[v_1(t, d_1 -1, \bs+\bbe_1)] \\ &\ge \ex_{1, t-1}[v_1(t, d_1, \bs -\bbe_j)] - \ex_{1, t-1}[v_1(t, d_1-1, \bs-\bbe_j+\bbe_1)].
\end{align*}
On the event $\Pi_1$, we know that, by induction assumption, we only need to consider the case that the seller one accepts the demand when the inventory level is at $d_1$, but not when it is at $d_1-1$. In this case, the LHS becomes $p$. For the RHS, consider the two cases that seller one accepts in both cases and only accepts when the inventory is $d_1$ but not $d_1-1$. In the first case, it becomes
\begin{align*}
\ex_{1, t}[v_1(t+1,  d_1-1, \bs -\bbe_j+\bbe_1)] - \ex_{1, t}[v_1(t+1, d_1-2, \bs -\bbe_j+\bbe_1)].
\end{align*}
Then the inequality follows from the condition that the seller accepts when the inventory and history is $( d_1-1, \bs -\bbe_j)$. In the second case, both the LHS and RHS become $p$. Now for the event $\Pi_j$, again, the one non-trivial case is similar. Hence, the LHS becomes,
\begin{align*}
&\ex_{1, t}[ v_1(t+1,  d_1, \bs+\bbe_j)] - \ex_{1, t}[v_1(t+1,  d_1-1, \bs + \bbe_j+\bbe_1)]\\ & \ge \ex_{1, t}[v_1(t+1,  d_1, \bs )] - \ex_{1, t}[v_1(t+1,  d_1-1, \bs+\bbe_1 )],
\end{align*}
and the inequality thus follows by induction.

This concludes the proof. $\Box$

\end{appendix}

\bibliographystyle{abbrv}
\bibliography{Lu}

\begin{thebibliography}{10}

\bibitem{BitranMondschein95}
G.~R. Bitran and S.~V. Mondschein.
\newblock An application of yield management to the hotel industry considering
  multiple day stays.
\newblock {\em Operations Research}, 43(3):427--443, 1995.

\bibitem{DenBoer2015}
A.~V. Den~Boer.
\newblock Dynamic pricing and learning: Historical origins, current research,
  and new directions.
\newblock {\em Surveys in Operations Research and Management Science},
  20(1):1--18, 2015.

\bibitem{FengGallego1995}
Y.~Feng and G.~Gallego.
\newblock Optimal starting times for end-of-season sales and optimal stopping
  times for promotional fares.
\newblock {\em Manage. Sci.}, 41(8):1371--1391, Aug. 1995.

\bibitem{FengGallego2000}
Y.~Feng and G.~Gallego.
\newblock Perishable asset revenue management with markovian time dependent
  demand intensities.
\newblock {\em Management Science}, 46(7):941--956, 2000.

\bibitem{FengXiao2000}
Y.~Feng and B.~Xiao.
\newblock Optimal policies of yield management with multiple predetermined
  prices.
\newblock {\em Operations Research}, 48(2):332--343, 2000.

\bibitem{doi:10.1287/mnsc.2013.1821}
G.~Gallego and M.~Hu.
\newblock Dynamic pricing of perishable assets under competition.
\newblock {\em Management Science}, 60(5):1241--1259, 2014.

\bibitem{GallegovanRyzin1994}
G.~Gallego and G.~van Ryzin.
\newblock Optimal dynamic pricing of inventories with stochastic demand over
  finite horizons.
\newblock {\em Manage. Sci.}, 40(8):999--1020, Aug. 1994.

\bibitem{LinLuYao2008}
G.~Y. Lin, Y.~Lu, and D.~D. Yao.
\newblock The stochastic knapsack revisited: Switch-over policies and dynamic
  pricing.
\newblock {\em Oper. Res.}, 56(4):945--957, July 2008.

\bibitem{osborne1994course}
M.~J. Osborne and A.~Rubinstein.
\newblock {\em A course in game theory}.
\newblock The MIT press, 1994.

\bibitem{PapastavrouRajagopalanKleywegt1996}
J.~D. Papastavrou, S.~Rajagopalan, and A.~J. Kleywegt.
\newblock The dynamic and stochastic knapsack problem with deadlines.
\newblock {\em Management Science}, 42(12):1706--1718, 1996.

\bibitem{VanSlykeYoung2000}
R.~Van~Slyke and Y.~Young.
\newblock Finite horizon stochastic knapsacks with applications to yield
  management.
\newblock {\em Operations Research}, 48(1):155--172, 2000.

\bibitem{ZhaoZheng2000}
W.~Zhao and Y.-S. Zheng.
\newblock Optimal dynamic pricing for perishable assets with nonhomogeneous
  demand.
\newblock {\em Management Science}, 46(3):375--388, 2000.

\end{thebibliography}

\end{document}